\documentclass[11pt]{amsart}
\usepackage[english]{babel}
\usepackage{graphicx}
\usepackage{framed}
\usepackage[normalem]{ulem}
\usepackage{amsmath}
\usepackage{amsthm}
\usepackage[mathscr]{euscript}

\usepackage{amssymb}
\usepackage{amsfonts}
\usepackage{enumerate}
\usepackage[utf8]{inputenc}
\usepackage[top=1 in,bottom=1in, left=1 in, right=1 in]{geometry}

\newcommand{\R}{\mathbb{R}}
\newcommand{\N}{\mathbb{N}}

\theoremstyle{definition}
\newtheorem{theorem}{Theorem}[section]
\newtheorem{corollary}[theorem]{Corollary}
\newtheorem{lemma}[theorem]{Lemma}
\newtheorem{proposition}[theorem]{Proposition}
\theoremstyle{definition}
\newtheorem{definition}[theorem]{Definition}
\newtheorem{remark}[theorem]{Remark}

\title[Finitary Codings]{Finitary Codings Of Finite Expectation Between Shifts Over Free Groups}
\author{James O'Quinn}
\date{}

\begin{document}

\maketitle
Some of the most central and well studied examples of  measurable dynamical systems are Bernoulli shifts, the theory of which is closely tied to dynamical entropy. The Ornstein isomorphism theorem states that Bernoulli shifts over the integers are classified, up to conjugacy, by their entropy. Since Ornstein's original proof of this theorem, which appeared in 1970, a great amount of interest and effort has been invested in extending it in various ways. \par
Firstly, one can attempt to show that shifts over a wider class of groups can be classified. Much progress has been made in this direction. Indeed, Ornstein and Weiss \cite{OW} generalized Ornstein's theorem to Bernoulli shifts over amenable groups. The breakthrough of a theory of entropy for sofic groups by Bowen \cite{BO} with improvements from Kerr and Li \cite{KL} allows one to distinguish shifts over sofic groups. Combining this with a result of Seward \cite{SE}, which builds on work of Bowen \cite{BO2}, we have that Bernoulli shifts over sofic groups are classified by their base entropy. \par 
Another way to extend the Ornstein isomophism theorem is to replace classification up to conjugacy by an even stronger form of equivalence. One such strengthening is finitary equivalence for Bernoulli shifts, which is defined by the existence of maps called finitary codings.  Roughly speaking, two Bernoulli shifts are finitarily equivalent if there is a conjugacy map between them, called a finitary coding isomorphism, such that the image of almost every point under this map is determined by an algorithm that is implemented by finite blocks of code at each step.  Coding information for a map between shifts is essentially topological: basic open sets in the product topology are the sets of points which share a given finite block of code. This means that some of the tools and techniques used in the theory of topological dynamics can be applied to the study of finitary codings. One key result is the theorem of Keane and Smorodinsky \cite{KS} stating that Bernoulli shifts over the integers with equal entropy are finitarily equivalent, strengthening the Ornstein isomorphism theorem.\par
The question is: how far can the notion of finitary equivalence be refined so as to still hold for Bernoulli shifts of equal entropy? Schmidt \cite{SC} showed, following the work of Parry and Schmidt \cite{PS}, that the additional assumption of having finite expected code length implies that the shifts must have equivalent distributions. Let $\mathbb{Z}\curvearrowright(X_{p},\mu_{p})$ and $\mathbb{Z}\curvearrowright(X_{q},\mu_{q})$ be Bernoulli shifts with initial distributions $p$ and $q$ respectively, and let $\phi:(X_{p},\mu_{p})\to(X_{q},\mu_{q})$ be a finitary coding (a measure preserving conjugacy which is also a continuous function a.e.\ on $X_{p}$). Thus, for a.e.\ $x\in X_{p}$, there are positive integers $a$ and $m$ such that, if $x'\in X_{p}$ satisfies $x_{i}=x'_{i}$ for all $-m\leq i\leq a$, we have that $\phi(x)_{0}=\phi(x')_{0}$. With $a(x)$ and $m(x)$ denoting the smallest positive integers that guarantee the above property for a given $x\in X_{p}$, we say that $\phi$ has finite expected code length if $\int \bigg(a(x)+b(x)\bigg)d\mu_{p}(x)<\infty$. One of the main technical consequences of $\phi$ having finite expected code length is the following observation due to Krieger in \cite{Kr}: for a.e. $x\in X_{p}$, we can find a $m\in\mathbb{N}$ (respectively $a\in \mathbb{N}$) such that, if $x'\in X_{p}$ is such that $x_{i}=x'_{i}$ for $i\in [-m,\infty)$ (respectively $i\in(\infty,a]$), we have $\phi(x)_{i}=\phi(x)_{i}$ for all $i\in [0,\infty)$ (respectively $i\in(-\infty,0)$). Utilizing the information cocycle, a perturbative extension of the classical information function introduced in \cite{BS}, Parry and Schmidt proved in \cite{PS} that there is a positive measure subset of $X_{p}$ on which the relative information created by a certain group generated by finitely supported permutations of $X_{p}$ is unchanged by the finitary isomorphism. Schmidt showed that this property can be extended to all of $X_{p}$, which is the key to equality of $p$ and $q$ up to some permutation. \par

Recently, Seward \cite{SE} generalized the result of Keane and Smorodinsky to all countably infinite groups, renewing some interest in the theory of finitary codings. In this paper, we aim to generalize Schmidt's result about finitary codings with finite expected code length to free groups. While amenable groups may seem, at first glance, to be the most natural setting for applying the techniques in \cite{SC} due to the asymptotic averaging used in the definition of the beta function, it is actually freeness which underpins the techniques described in the previous paragraph. Freeness allows for the implementation of a Turing-machine-like scheme inside of the group: we can isolate an infinite set inside the group such that translation only moves one element outside the set. Thus, the information cocycle tells us that shifting a point in the space only changes the relative information with respect to this set by a single symbol. This in turn allows us to find a formula for the beta function, similar to the one in \cite{SC}, as an asymptotic average of the information cocycle over the iterates of a shift. Also, Krieger's proof of the observation from \cite{Kr} relies on the pointwise ergodic theorem, which only applies to amenable groups. In Proposition 1.3 below an alternate proof of Kreiger's observation will be given without the use of the pointwise ergodic theorem, or any other tools related to amenability. \par 
\emph{Acknowledgements}. The author greatly appreciates the helpful comments and advice from David Kerr. The research was partially funded by the Deutsche Forschungsgemeinschaft (DFG, German Research Foundation) under Germany’s Excellence Strategy – EXC
2044 – 390685587, Mathematics M¨unster – Dynamics – Geometry – Structure; the Deutsche
Forschungsgemeinschaft (DFG, German Research Foundation) – Project-ID 427320536 – SFB
1442, and ERC Advanced Grant 834267 - AMAREC

\section{Finitary Codings}
Let $F_{\ell}$ be the free group on $\ell$ generators, and write $S$ for this set of generators. We define the word length function $|\cdot|$ on $F_{\ell}$ as follows: given some word $g\in F_{\ell}\setminus\{e\}$, set $|g|=\min\{m\in \mathbb{N}\mid\exists g_{1},\dots,g_{m}\in S\cup S^{-1} \text{ such that }g=g_{1}\dots g_{m}\}$ and $|e|=0$. Let $r\in \mathbb{N}$ and write $B(r)$ for the set of all words in $F_{\ell}$ whose word length is at most $r$. Define $p=[P_{1},\dots,P_{m}]$ to be a vector whose entries are strictly positive real numbers that sum up to $1$. We get a probability space on the set $\{1,\dots,m\}$ using the probability measure $P(i)=P_{i}$ for every $i\in \{1,\dots,m\}$. This gives another probability space $X_{p}=\{1,\dots,m\}^{F_{\ell}}$ equipped with the product probability measure $\mu_{p}=\prod_{g\in F_{\ell}}P^{g}$ when setting $P_{g}=P$ for all $g\in F_{\ell}$. We will use the following notation to denote basic open (cylinder) sets in $X_{p}$: given some $A\subset F_{\ell}$ and a set $\{i_{g}\}_{g\in A}$ where $i_{g}\in \{1,\dots,m\}$ for all $g\in A$, we set $[i_{g}\mid g\in A]=\{x\in X_{p}\mid x_{h}=i_{h} \text{ for every } h\in A\}$. We get a probability measure preserving action
$F_{\ell}  \curvearrowright (X_{p},\mu_{p})$ by setting $(g x)_{h}=x_{g^{-1}h}$ for all $x\in X_{p}$ and all $g,h\in F_{\ell}$.  Set $\alpha_{p}=\{[i_{g}\mid g\in\{e\}]: i_{g}\in\{1,\dots,m\}\}$ where e is the empty word (identity element) in $F_{\ell}$. We call $\alpha_{p}$ the identity partition of $X_{p}$. For the rest of the paper, fix an $a\in S$ and set $W_{a}=\{wa\in F_{\ell}:|wa|=|w|+1\}\cup\{e\}$. We can see that $W_{a}$ is the set of reduced words that end with an $a$. We define the \emph{past algebra} $\mathscr{A}_{p,a}$ of $F_{\ell}  \curvearrowright (X_{p},\mu_{p})$ over $W_{a}$ to be the sub-$\sigma$-algebra generated by $\{g^{-1}\alpha_{p}\mid g\in W_{a}\}$. 
\par
We can also view $X_{p}$ as a topological space. Indeed, we can equip $X_{p}$ with the topology generated by the sub-basis $\{[i_{g}\mid g\in\{h\}\mid h\in F_{\ell}\text{ and } i_{g}\in \{1,\dots,m\}\}$, under which it is a Cantor set (unless $p$ is trivial).
\par
Let $q=[Q_{1},\dots,Q_{n}]$ be a second probability vector. This induces a Bernoulli  shift $F_{\ell} \curvearrowright (X_{q},\mu_{q})$ with identity partition $\alpha_{q}$. Let $\mathscr{A}_{q,a}$ be the sub-$\sigma$-algebra generated by $\{g^{-1}\alpha_{q}\mid g\in W_{a}\}$. We also equip $X_{q}$ with the topology generated by the sub-basis $\{[i_{g}\mid g\in\{h\}]\mid h\in F_{\ell}\text{ and } i_{g}\in \{1,\dots,n\}\}$.

\begin{definition}A measurable map $\phi :(X_{p},\mu_{p})\to (X_{q},\mu_{q})$ is called a \textbf{finitary coding} if there exists a null set $E\subset X_{p}$ such that the restriction of $\phi$ to the set $X_{p}\setminus E$ is a continuous map. Thus, for almost every $x\in X_{p} $, we have that $\phi(x)_{e}$ is determined by a ball $B(r)\subset F_{\ell}$ in the following way: there exists an $r\in\mathbb{N}$ depending on $x$ such that if $x'\in X_{p}$ such that $x'_{g}=x_{g}$ for every $g\in B(r)$, then $\phi(x)_{e}=\phi(x')_{e}$. Let  $r_{\phi}(x)$ be the minimum radius $r$ for which $B(r)$ determines $\phi(x)_{e}$. We call $\phi$ a \textbf{finitary isomorphism} if it is a conjugacy and a measure preserving bijection and $\phi$ and $\phi^{-1}$ are finitary codings. Finally, a finitary isomorphism  $\phi$ has \textbf{finite expectation} if, setting $v_{p}(x)=|B(r_{\phi}(x))|$, we have $\int v_{\phi}(x)d\mu_{p}<\infty$ and $\int v_{\phi^{-1}}(x)d\mu_{q}<\infty$.

\end{definition}

\begin{definition} Given a finitary coding $\phi:X_{p}\to X_{q}$, we can define functions $m_{\phi}:X_{p}\to \mathbb{Z}$ and $a_{\phi}:X_{p}\to \mathbb{Z}$ in the following way:
\begin{align*} m_{\phi}(x)= \underset{g\in W_{a}}{\sup }(r_{\phi}(g^{-1}x)-|g|)\text{ and }
\end{align*}
\begin{align*} a_{\phi}(x)= \underset{g\in F_{\ell}\setminus W_{a}}{\sup }(r_{\phi}(g^{-1}x)-|g|). \end{align*}
\end{definition}

The functions $m_{\phi}$ and $a_{\phi}$ were first introduced in \cite{PS} as the functions $m^{*}$ and $a^{*}$ following an idea of \cite{Kr}. The idea is that, for some point $x\in X_{p}$, $m_{\phi}$ (respectively $a_{\phi}$) determines the smallest radius of a ball centered at the identity in $F_{\ell}$ that is needed to determine $\phi(x)_{g}$ for all $g\in W_{a}$ (respectively $g\in F_{\ell}\setminus W_{a}$). The next proposition shows that this property holds for a.e. $x\in X_{p}$.

\begin{proposition} If $\phi:X_{p}\to X_{q}$ is a finitary isomorphism with finite expectation, then the functions $m_{\phi}$ and $a_{\phi}$ are finite a.e.\par
    \begin{proof}

 First, it will be shown that $m_{\phi}$ is finite a.e. Let $A=\bigcap_{i=0}^{\infty}\bigcup_{g\in W_{a}}\{x\mid r_{\phi}(g^{-1}x)-|g|>i\}$. It suffices to show that $\mu_{p}(A)=0$. Let $\epsilon>0$ be given. The next part of the proof will proceed in two cases: when $\ell=1$ and when $\ell>1$.\par
$\textbf{Case 1}$: $\ell=1$.\par
Let $a$ be the generator of $F_{\ell}$. Since $\phi$ has finite expectation, we have that $\sum_{n=1}^{\infty}\mu_{p}(\{x\mid v_{\phi}(x)>n\})=\int v_{\phi}(x) d\mu_{p} <\infty$. Since $F_{\ell}\curvearrowright X_{p}$ is measure preserving, we have that $\mu_{p}(\{x\mid v_{\phi}(a^{-n}x)>n\})=\mu_{p}(\{x\mid v_{\phi}(x)>n\})$ for all $n\in \mathbb{N}$. So, there exists an $N\in\N$ such that $N>3$ and 
    \begin{align*}
    \sum_{n=N}^{\infty}\mu_{p}(\{x\mid v_{\phi}(a^{-n}x)>n\})=\sum_{n=N}^{\infty}\mu_{p}(\{x\mid v_{\phi}(x)>n\})<\epsilon.
\end{align*}
\par Now, it will be shown that $A\subset \bigcup_{n=N}^{\infty}\{x\mid v_{\phi}(a^{-n}x)>n\}$. Suppose by way of contradiction that there exists a $y\in A$ such that $y\notin \bigcup_{n=N}^{\infty}\{x\mid v_{\phi}(a^{-n}x)>n\}$. Since $y\in A$, there exists an integer $k>N$ such that $r_{\phi}(a^{-k}x)-k>N$. Since $y\notin \bigcup_{n=N}^{\infty}\{x\mid v_{\phi}(a^{-n}x)>n\}$, we have that $v_{\phi}(a^{-k}y)\leq k$. However, $v_{\phi}(x)=2r_{\phi}(x)-1$ which implies that $2(k+N)-1<v_{\phi}(a^{-k}y)\leq k$, a contradiction. Thus $A\subset \bigcup_{n=N}^{\infty}\{x\mid v_{\phi}(g_{n}^{-1}x)>n\}$. Now we will continue with the second case.
\par
$\textbf{Case 2}$: $\ell>1$. \par
   Let $\{g_{n}\}_{n=1}^{\infty}$ be an enumeration of $W_{a}$ such that $|g_{n}|\leq |g_{n+1}|$ for all $n\in \mathbb{N}$. For an example of such an enumeration, take the enumeration produced by arbitrarily ordering words of a given word length and then enumerating these orderings sequentially across all word lengths. In fact, these are the only such enumerations. This produces an enumeration with the property $n\leq (2\ell-1)^{|g_{n}|+1}$ for all $n\in\mathbb{N}$. Indeed, we can see that $|g_{n}|\leq |g_{n+1}|$ for all $n\in \mathbb{N}$ implies that the enumeration must also preserve word length. Thus the enumeration of a given word $g$ is bounded by the number of all possible words in $W_{a}$ with word length at most $|g|$, which is bounded by $(2\ell-1)^{|g|+1}$.  Since $\phi$ has finite expectation, we have that $\sum_{n=1}^{\infty}\mu_{p}(\{x\mid v_{\phi}(x)>n\})=\int v_{\phi}(x) d\mu_{p} <\infty$. Since $F_{\ell}\curvearrowright X_{p}$ is measure preserving, we have that $\mu_{p}(\{x\mid v_{\phi}(g_{n}^{-1}x)>n\})=\mu_{p}(\{x\mid v_{\phi}(x)>n\})$. So there exists an $N\in\N$ such that $N>3$ and 
    \begin{align*}
    \sum_{n=N}^{\infty}\mu_{p}(\{x\mid v_{\phi}(g_{n}^{-1}x)>n\})=\sum_{n=N}^{\infty}\mu_{p}(\{x\mid v_{\phi}(x)>n\})<\epsilon.
\end{align*}
\par Now, it will be shown that $A\subset \bigcup_{n=N}^{\infty}\{x\mid v_{\phi}(g_{n}^{-1}x)>n\}$. Suppose by way of contradiction that there exists a $y\in A$ such that $y\notin \bigcup_{n=N}^{\infty}\{x\mid v_{\phi}(g_{n}^{-1}x)>n\}$. Since $k\leq (2\ell-1)^{|g_{k}|+1}$  for all $k\in\mathbb{N}$ and $y\in A$, we can find a $k$ to make $|g_{k}|$, and in turn $r_{\phi}(g_{k}^{-1}y)$, as large as we want. Thus, we can find a large enough $k$ to make $\log_{2\ell-1}((2\ell-1)^{r_{\phi}(g_{k}^{-1}y)}-1)$ as close to $r_{\phi}(g_{k}^{-1}y)$ as we need.
     So, since $y\in A$, we can find a large enough integer $k>N$ such that, 
\begin{enumerate}
\item $r_{\phi}(g_{k}^{-1}y)-|g_{k}|>N$
  
\item $r_{\phi}(g_{k}^{-1}y)-\log_{2\ell-1}((2\ell-1)^{r_{\phi}(g_{k}^{-1}y)}-1)\leq\log_{2\ell-1}(\frac{2\ell}{2\ell-2})$. 
    
\end{enumerate}
Note that $v_{\phi}(y)=1+\frac{2\ell((2\ell-1)^{r_{\phi}(y)}-1)}{2\ell-2}$ by using the usual formula for the cardinaility of a ball in $F_{\ell}$. By hypothesis, $v_{\phi}(g_{k}^{-1}y)\leq k$, which implies that 
\begin{align*}
N+|g_{k}|&\stackrel{(1)}{\leq} r_{\phi}(g_{k}^{-1}y) \\ 
    &\stackrel{(2)}{\leq} \log_{2\ell-1}((2\ell-1)^{r_{\phi}(g_{k}^{-1}y)}-1) +\log_{2\ell-1}(\frac{2\ell}{2\ell-2})\\
    &\leq\log_{2\ell-1}(v_{\phi}(g_{k}^{-1}y))\\
    &\leq \log_{2\ell-1}(k) \,\,\,\,\,\,\,\,\,\,\,\,\,\,\,\,\,\,\,\,\,\,\,\text{  since }y\notin \bigcup_{n=N}^{\infty}\{x\mid v_{\phi}(g_{n}^{-1}x)>n\} \\
    & \leq |g_{k}|+1 
\end{align*} which gives a contradiction. Thus, $A\subset \bigcup_{n=N}^{\infty}\{x\mid v_{\phi}(g_{n}^{-1}x)>n\}$.\par
So, in both cases, we have the following:
 \begin{align*}
    \mu_{p}(A)\leq \mu_{p}\bigg(\bigcup_{n=N}^{\infty}\{x\mid v_{\phi}(g_{n}^{-1}x)>n\}\bigg)\leq \sum_{n=N}^{\infty}\mu_{p}(\{x\mid v_{\phi}(g_{n}^{-1}x)>n\})<\epsilon.
\end{align*} 
This proves that $m^{*}$ is finite almost everywhere. The proof that $a^{*}$ 
is finite almost everywhere follows in a similar manner with a few modifications. For the case $\ell=1$, just replace $W_{a}$ with $F_{\ell}\setminus W_{a}$ and the proof proceeds in the same way as above. For the case $\ell>1$, we can enumerate $F_{\ell}\setminus W_{a}$ using the same method as before, only this time we get the property that $n\leq 3(2\ell-1)^{|g_{n}|+1}$ for all $n\in\mathbb{N}$. By taking the same $N$ and $k$ as given above, we get $N+|g_{k}|\leq |g_{k}|+1+ \log_{2\ell-1}3$ which is still a contradiction.
 \end{proof}
\end{proposition}
The main technical uses of Proposition 1.3 are detailed in the following corollary.
\begin{corollary}
If $\phi:X_{p}\to X_{q}$ is a finitary coding with finite expectation, then the following hold:
\begin{enumerate}
    \item For a.e. $x\in X_{p}$, there exists a ball $B\subset F_{\ell}$ centered at the identity element such that if $x'\in X_{p}$ has the property that $x'_{g}=x_{g}$ for all $g\in B\cup W_{a}$, then $\phi(x')_{g}=\phi(x)_{g}$ for all $g\in W_{a}$.
    \item For a.e. $x\in X_{p}$, there exists a ball $B\subset F_{\ell}$ centered at the identity element such that if $x'\in X_{p}$ has the property that $x'_{g}=x_{g}$ for all $g\in B\cup F_{\ell}\setminus W_{a}$, then $\phi(x')_{g}=\phi(x)_{g}$ for all $g\in F_{\ell}\setminus W_{a}$.
\end{enumerate}
\end{corollary}
\section{Conditional Information Functions}
In this section, we will define the conditional information functions that we obtain by factoring a measure space via a sub-$\sigma$-algebra using measurable partitions. For this section, let $(X,\mathscr{F},\mu)$ be a probability space. The following definitions and facts can be found in \cite{ME} and \cite{BS}.

\begin{definition}
A collection $\alpha$ of disjoint measurable subsets of $X$ is a \textbf{measurable partition} if $\bigsqcup_{A\in\alpha}A=X$ and has the following separation property:
\begin{itemize}
    \item there is a countable collection $\{B_{n}\}_{n\in\mathbb{N}}$ of measurable sets, each of which is a union of sets from $\alpha$, such that if distinct $A_{1},A_{2}\in \alpha$ then there exists an $n\in \mathbb{N}$ such that either $A_{1}\subset B_{n}$ and $A_{2}\subset X\setminus B_{n}$ or vice versa.
\end{itemize}
Note that the measurable partitions defined above need not be countable. Also, we will call the set in the partition that contains a particular $x\in X$ the \textbf{atom} of $x$ in $\alpha$.
\end{definition}
\begin{definition}
A measurable partition $\alpha$ induces a $\textbf{factor space}$ $(X_{\alpha},\mathscr{F}_{\alpha},\mu_{\alpha})$ of $(X,\mathscr{F},\mu)$ where $X_{\alpha}=\alpha$, $\mathscr{F}_{\alpha}$ is the $\sigma$-algebra generated by all $\mathscr{F}$-measurable sets which are a union of sets from $\alpha$, and $\mu_{\alpha}$ is the restriction of $\mu$ to $\mathscr{F}_{\alpha}$. This also induces a $\textbf{factor projection}$ $N_{\alpha}:X\to X_{\alpha}$ which sends a point in $X$ to the set in $\alpha$ which contains it. 
\end{definition}
The following theorem appears as Theorem 1.14 in \cite{ME}.
\begin{theorem}
Let $\alpha$ be a measurable partition of the probability space $(X,\mathscr{F},\mu)$, and let $N_{\alpha}$ be the associated projection for the factor space $(X_{\alpha},\mathscr{F}_{\alpha},\mu_{\alpha})$. For $\mu_{\alpha}$-a.e. $A\in\alpha$ there exists a probability measure $\mu^{\alpha}(\cdot|A)$ on $\mathscr{F}$ such that:
\begin{enumerate}
    \item $\mu^{\alpha}(X\setminus N_{\alpha}^{-1}(A) |A)=0$.
    \item For any $F\in\mathscr{F}$, $A\to \mu^{\alpha}(F|A)$ is an $\mathscr{F}_{\alpha}$-measurable function defined $\mu_{\alpha}$ almost everywhere on $X_{\alpha}$.
    \item For any $F\in\mathscr{F}$, we have  \begin{align*}
        \mu(F)&=\int_{X_{\alpha}}\mu^{\alpha}(F|A)d\mu_{\alpha}(A)\\
        &=\int_{X}\mu^{\alpha}(F|N_{\alpha}(x))d\mu(x).
    \end{align*}
\end{enumerate}
\end{theorem}
The main use we will have for these conditional measures in this paper is that we can factor a measure space over a sub-$\sigma$-algebra. We can do this by obtaining a measurable partition from a given sub-$\sigma$-algebra. We will detail this procedure below, but first we need a definition and a couple of theorems.
\begin{definition}
Let $(X,\mathscr{F},\mu)$ be a probability space. We say a countable collection of measurable sets $\Gamma$ is a \textbf{basis} for $(X,\mathscr{F},\mu)$ if it satisfies the following two conditions:
\begin{enumerate}
    \item For any two distinct points $x_{1},x_{2}\in X$, there exists a set $B\in\Gamma$  such that either $x_{1}\in B$ and $x_{2}\notin B$ or vice versa.
    \item The completion of the $\sigma$-algebra generated by $\Gamma$ is $\mathscr{F}$.
\end{enumerate}
\end{definition}

The following theorems appear in \cite{ME}.
\begin{theorem}(Rokhlin) If $\alpha$ is a measurable partition of a probability space $(X,\mathscr{F},\mu)$ with a basis, then the factor space $(X_{\alpha},\mathscr{F}_{\alpha},\mu_{\alpha})$ also has a basis.

\end{theorem}
\begin{theorem}
If $(X,\mathscr{F},\mu)$ is a probability space with a basis and $\mathscr{A}$ is a sub-$\sigma$-algebra of $\mathscr{F}$, then $(X,\mathscr{A},\mu)$ also has a basis.
\end{theorem}
\begin{remark}

Given a sub-$\sigma$-algebra $\mathscr{A}$ of $\mathscr{F}$, we can get a measurable partition $\alpha$ of $X$ from $\mathscr{A}$ in the following way: let $\{A_{n}\}_{n\in \mathbb{N}}$ be a countable basis for $\mathscr{A}$ and define $\alpha$ to the partition consisting of the sets of the form $\bigcap_{n\in\mathbb{N}}\hat{A}_{n}$ where either $\hat{A}_{n}=A_{n}$ or $\hat{A}_{n}=X\setminus A_{n}$. This gives a unique measurable partition on $X$ which is independent of the choice of $\{A_{n}\}_{n\in \mathbb{N}}$ up to a null set in $X$, i.e., given two measurable partitions $\alpha$ and $\beta$ consisting of sets in $\mathscr{A}$, there exists a null set $N$ such that $\{A\setminus N\mid A\in\alpha\}=\beta$. The atom that contains a given $x\in X$ in this partition will be denoted by $[x]_{\mathscr{A}}$. Also, we will use $\mu^{\mathscr{A}}$ to denote the family of conditional measures $\{\mu^{\alpha}(\cdot|A)\}_{A\in\alpha}$ of $\mu$ induced by $\mathscr{A}$ and $N_{\mathscr{A}}$ to denote the associated projection.
\end{remark}

\begin{definition}
Let two sub-$\sigma$-algebras $\mathscr{A}$ and $\mathscr{B}$ of $\mathscr{F}$ be given. We define the $\textbf{conditional}$ $\textbf{information}$ of $\mathscr{B}$ given $\mathscr{A}$ and with respect to $\mu$ to be the function
\begin{align*}
    I_{\mu}(\mathscr{B}|\mathscr{A})(x)=-\log \mu^{\mathscr{A}}([x]_{\mathscr{B}}\mid [x]_{\mathscr{A}}), \text{ }x\in X.
\end{align*}
\end{definition}

We can also define a conditional expectation as well.

\begin{definition}
Given a sub-$\sigma$-algebra $\mathscr{A}$ of $(X,\mathscr{F},\mu)$, We can define the \textbf{conditional expectation} of $\mathscr{A}$ of $f\in L^{1}(X,\mathscr{F},\mu)$ to be the function \begin{align*}
E(f|\mathscr{A})(x)=\int fd\mu^{\mathscr{A}}(\cdot|[x]_{\mathscr{A}})
\end{align*}
defined $\mu$-a.e.
\end{definition}
\section{The Information Cocycle}
In this section we will define the information cocycle, which is the main tool we will use to pass conditional information from the past algebra of one shift to the past algebra of the other under the finitary coding. We return to using the notation and conventions from Section 1 and Section 2. 
\begin{definition}
We say an automorphism $V:X_{p}\to X_{p}$ is \textbf{locally finite} if $|\{g\in F_{\ell}\mid V(x)_{g}\neq x_{g}\}|<\infty$ for a.e. $x\in X_{p}$. Let $L_{p}$ denote the group of all locally finite automorphisms of $(X_{p},\mu_{p})$. Let $G_{p}$ be the group generated by $L_{p}$ and $F_{\ell}$. Let $a$ be a generator of $F_{\ell}$ and let $L_{p,a}$ be the subgroup of $L_{p}$ such that $V(x)_{g}=x_{g}$ whenever $g\in W_{a}\setminus \{a^{k}\mid k\in\mathbb{Z}\}$ for a.e. $x\in X_{p}$ when $V\in L_{p,a}$. The \textbf{information cocycle} $J$ conditioned by a sub-$\sigma$-algebra $\mathscr{A}$ of $(X_{p},\mu_{p})$ is the function
\begin{align*} J(\mathscr{A},V)=I_{\mu_{p}}(\mathscr{A}|V^{-1}\mathscr{A})-I_{\mu_{p}}(V^{-1}\mathscr{A}|\mathscr{A})-\log E\bigg{(}\frac{d\mu_{p}V^{-1}}{d\mu_{p}}|\mathscr{A}\bigg{)}\circ g
\end{align*}
where $V\in G_{p}$, $I$ is the conditional information function, and $E$ is the conditional expectation, as given in \cite{BS}.
It is not too hard to see that $J$ satisfies the following cocycle equation:
\begin{align*}
    J(\mathscr{A},VW)= J(\mathscr{A},V)\circ W+ J(\mathscr{A},W)
\end{align*}
a.e. on $X_{p}$.

\end{definition}
The next proposition gives some useful formulas for calculating the informaiton cocycle.
\begin{proposition}
Let $n\in \mathbb{N}$ and $a$ be a standard generator of $F_{\ell}$. Then, we have the following formula for the information cocycle conditioned by the past algebra $\mathscr{A}_{p,a}$: 
\begin{align}\label{eqn1}
J(\mathscr{A}_{p,a},a^{n})(x)=\sum_{i=0}^{n-1}I_{\mu_{p}}(\mathscr{A}_{p,a}|a^{-1}\mathscr{A}_{p,a})(a^{i}x)=\sum_{i=0}^{n-1}-\log(P_{x_{a^{i}}})
\end{align} for $\mu_{p}$-a.e. $x\in X_{p}$. Moreover, given $V\in L_{p,a}$, we have for $\mu_{p}$-a.e. $x\in X_{p}$ that
\begin{align}\label{eqn2}
J(\mathscr{A}_{p,a},V)(x)=\log\bigg{(}\prod_{g\in W_{a}}\frac{P_{(Vx)_{g}}}{P_{x_{g}}}\bigg{)}.
\end{align}
\end{proposition}
\begin{proof}
 To produce equation (1), it only needs to be shown that $J(\mathscr{A}_{p,a},a)(x)=-\log(P_{x_{e}})$ since the cocycle identity can be applied recursively with this equation to produce the formula for $J(\mathscr{A}_{p,a},a^{n})(x)$, $n>1$. Since shifting by $a$ is a measure preserving transformation, we have that \begin{align*}
     J(\mathscr{A},a)=I_{\mu_{p}}(\mathscr{A}|a^{-1}\mathscr{A})-I_{\mu_{p}}(a^{-1}\mathscr{A}|\mathscr{A}).
 \end{align*} It is not too hard to see that $a^{-1}\mathscr{A}$ is the sub-$\sigma$-algebra of $\mathscr{A}$ generated by the collection $\{g^{-1}\alpha_{p}\mid g\in W_{a}\}$. Let $x\in X_{p}$. We have that $[x]_{a^{-1}\mathscr{A}}=\{y\in X_{p}\mid y_{g}=x_{g}\text{ }\forall g\in W_{a}\}$ and $[x]_{\mathscr{A}}=\{y\in X_{p}\mid y_{g}=x_{g}\text{ }\forall g\in W_{a}\cup \{e\}\}$. Thus, we have that $[x]_{\mathscr{A}}\subset [x]_{a^{-1}\mathscr{A}}$. This implies that $I_{\mu_{p}}(a^{-1}\mathscr{A}|\mathscr{A})=-\log \mu_{p}^{\mathscr{A}}([x]_{a^{-1}\mathscr{A}}\mid [x]_{\mathscr{A}})=0$. On the other hand, we can see that $I_{\mu_{p}}(\mathscr{A}|a^{-1}\mathscr{A})=-\log \mu_{p}^{\mathscr{a^{-1}A}}([x]_{\mathscr{A}}\mid [x]_{a^{-1}\mathscr{A}})=-\log(P_{x_{e}})$.\par
 
 For the proof of (2), let $V\in L_{p,a}$ and let $x\in X_{p}$ such that there is an $n\in \mathbb{N}$ such that $a^{n}Va^{-n}[x]_{\mathscr{A}_{a,p}}=[x]_{\mathscr{A}_{a,p}}$. Thus, we can apply the cocycle identity to $a^{n}Va^{-n}$ to produce the following equation:
\begin{align*}
   J(\mathscr{A}_{p,a},V)(x)=-J(\mathscr{A}_{p,a},a^{n})(Vx)-J(\mathscr{A}_{p,a},a^{-n})(a^{n}x).\tag{$*$}
\end{align*}

Applying the cocycle identity to $aa^{-1}=e$ and using equation (1), we can see that:
\begin{align*}
    J(\mathscr{A}_{p,a},a^{-1})(x)=\log(P_{x_{a^{-1}}})\tag{$**$}
\end{align*}
and
\begin{align*}
    J(\mathscr{A}_{p,a},a^{-n})(x)=\sum_{i=0}^{n-1} J(\mathscr{A}_{p,a},a^{-1})(a^{-i}x).  \tag{$***$}
\end{align*}

Hence,
\begin{align*}
  J(\mathscr{A}_{p,a},V)(x) &\stackrel{(*)}{=}-J(\mathscr{A}_{p,a},a^{n})(Vx)-J(\mathscr{A}_{p,a},a^{-n})(a^{n}x) \\
  &\stackrel{(***)}{=} -\sum_{i=0}^{n-1}I_{\mu_{p}}(\mathscr{A}_{p,a}|a^{-1}\mathscr{A}_{p,a})(a^{i}Vx)-\sum_{i=0}^{n-1} J(\mathscr{A}_{p,a},a^{-1})(a^{-i}x) \\
  &\stackrel{(**)}{=} \sum_{i=0}^{n-1}\log(P_{V(x)_{a^{i}}})-\sum_{i=0}^{n-1}\log(P_{x_{a^{n-i}}}) \\
  &= \sum_{i=0}^{n-1}(\log(P_{V(x)_{a^{i}}})-\log(P_{x_{a^{i}}}))\\
  &= \sum_{i=0}^{n-1}\log\bigg(\frac{P_{V(x)_{a^{i}}}}{P_{x_{a^{i}}}}\bigg)\\
  &=\log\bigg(\prod_{i=0}^{n-1}\frac{P_{V(x)_{a^{i}}}}{P_{x_{a^{i}}}}\bigg).
\end{align*}
So, we get that $J(\mathscr{A}_{p,a},V)(x)=\log\bigg{(}\prod_{g\in W_{a}}\frac{P_{(Vx)_{g}}}{P_{x_{g}}}\bigg)$ for $\mu_{p}$-a.e. $x\in X_{p}$.

\end{proof}
\begin{proposition}
There exists a measurable (co-boundary) map $f:X_{p}\to \mathbb{R}$ such that
\begin{align*}J(\mathscr{A}_{p,a},g)=J(\mathscr{A}_{q,a},\phi g\phi^{-1})\circ\phi+f\circ g-f \text{ a.e.}
\end{align*}
for every $g\in G_{p}$.
\end{proposition}

\begin{proof}
By Theorem 3.5 from \cite{BS}, it only needs to be shown that 
\begin{align*}
    I_{\mu_{p}}(\phi^{-1}\mathscr{A}_{q,a}|\mathscr{A}_{p,a})+I_{\mu_{p}}(\mathscr{A}_{p,a}|\phi^{-1}\mathscr{A}_{q,a})<\infty
\end{align*}
almost everywhere on $X_{p}$. Indeed, let S be the subset of X given by Proposition 1.3 such that $m_{\phi}(x)<\infty$ and $m_{\phi^{-1}}(x)<\infty$ and $m_{p}(S)=1$. First, it will be shown that $I(\phi^{-1}\mathscr{A}_{q,a}|\mathscr{A}_{p,a})<\infty$ on S. This can be shown by finding, for every $x \in S$, a common cylinder set in both $[x]_{\mathscr{A}_{p,a}}$ and $[x]_{\phi^{-1}\mathscr{A}_{q,a}}$. Indeed, since $m_{\phi}(x)<\infty$, there is a ball $B\subset F_{\ell}$ centered at the identity and a cylinder set $C=\{y\in X_{p}|y_{g}=x_{g} \text{ for all } g\in B\cup W_{a}\}$ such that, for every $y\in C$, $\phi(y)_{g}=\phi(x)_{g}$ for all $g\in W_{a}$. This implies that $C\subset [x]_{\mathscr{A}_{p,a}}$ and $C\subset \phi^{-1}[\phi (x)]_{\mathscr{A}_{q,a}}=[x]_{\phi^{-1}\mathscr{A}_{q,a}}$. Also, we have that $C$ has nonzero measure when $m_{p}$ is conditioned by $[x]_{\mathscr{A}_{p,a}}$. Thus, $I(\phi^{-1}\mathscr{A}_{q,a}|\mathscr{A}_{p,a})<\infty$. \par
Now, we will show $I(\mathscr{A}_{p,a}|\phi^{-1}\mathscr{A}_{q,a})<\infty$ for every $x\in S$. To show this, we can apply the same technique used in this first part of the proof, this time using the fact that $\phi^{-1}$ also has finite expectation. Indeed, given an $x\in S$, we can find a ball $B\subset F_{\ell}$ centered at the identity and a cylinder set $C=\{y\in X_{q}\mid y_{g}=\phi(x)_{g} \text{ for all } g\in B\cup W_{a}\}$ in $X_{q}$ such that $\phi^{-1}C\subset \phi^{-1}[\phi(x)]_{\mathscr{A}_{q,a}}=[x]_{\phi^{-1}\mathscr{A}_{q,a}}$ and $\phi^{-1}C\subset [x]_{\mathscr{A}_{p,a}}$. Using property (2.11) from \cite{BS}, we get
\begin{align*}
    \mu_{p}^{\phi^{-1}\mathscr{A}_{q,a}}(\phi^{-1}C|[x]_{\phi^{-1}\mathscr{A}_{q,a}})=(\phi^{-1}\mu_{p,a})^{\mathscr{A}_{q,a}}(C|[\phi(x)]_{\mathscr{A}_{q,a}})>0.
\end{align*}
The last inequality holds since $\phi$ is non-singular.
\end{proof}
\bigskip{}

\par

Now, we will find a positive measure set in $X_{p}$ on which the co-boundary factor $f$ given by Proposition 3.3 is constant a.e.
\par
Let $C$ be a cylinder set in $X_{p}$ over a ball $B_{N}$ centered at the identity element of radius N. Let $H_{C}^{+}=\{h\in L_{p}\mid  (hx)_{g}=x_{g} \text{ for all } g\in B_{N}\cup (F_{\ell}\setminus  W_{a}) \text{ and for a.e. }x\in X_{p}\}$ and $H_{C}^{-}=\{h\in L_{p}\mid  (hx)_{g}=x_{g} \text{ for all } g\in W_{a}\cup B_{N}  \text{ and for a.e. }x\in X_{p}\}$. Finally, let $H_{C}=H_{C}^{+}\cdot H_{C}^{-} $.\par

For the next proposition, we will use the following characterization of a weakly mixing action. For more on weak mixing, see [].

\begin{definition}
Let $G$ be a group, $(X,\mathscr{A},\mu)$ be a probability space, and $G\curvearrowright (X,\mathscr{A},\mu)$ be a probability measure preserving action. We say $G\curvearrowright X$ is \textbf{weakly mixing} if for every $\epsilon>0$ and for every finite $\Omega\subset \mathscr{A}$ there is a $g\in G$ such that 
\begin{align*}
     |\mu(gA\cap B)-\mu(A)\mu(B)| <\epsilon
\end{align*}
whenever $A,B\in \Omega$.
\end{definition}

\begin{proposition} Let $\mu_{C}$ be the probability measure on $C$ defined by $\mu_{C}(A)=\frac{\mu_{p}(A)}{\mu_{p}(C)}$ for all measurable $A\subset C$. The action of $H_{C}$ on $(C,\mu_{C})$ is weakly mixing.

\end{proposition}
\begin{proof} First, we will show that the condition in the definition of weak mixing holds on finite collections of cylinder sets in $C$. Then, we can approximate in measure every finite collection of measurable sets from $C$ by a finite collection of cylinder sets to get the condition in the definition of weak mixing to hold for all measurable sets in $C$.\par 

Let $\Omega=\{C_{1},C_{2},\dots, C_{n}\}$ be a finite collection of cylinder sets in $C$ such that for each $k\in \{1,2,\dots,n\}$, there is a finite set $A_{k}\subset F_{\ell}$ and $\{i_{g}^{k}\}_{ g\in A_{k}}\subset \{1,2,\dots,m\}$ such that $C_{k}=[i_{g}^{k}\mid g\in A_{k}]$. Let $B'$ be the smallest ball in $F_{\ell}$ centered at the identity such that for every $k\in\{1,2,\dots,n\}$ we have $A_{k}\subset B'$. Define $h^{+}\in H_{C}^{+}$ in the following way. First, for every $g\in (B'\setminus B_{N})\cap W_{a}$, we choose a $h_{g}\in W_{a}\setminus B'$ such that $h_{g}=h_{g'}$ only if $g=g'$ for all $g,g'\in (B'\setminus B_{N})\cap W_{a}$. Then, for $x\in C$, we set
\begin{align*}
h^{+}(x)_{f}=
    \begin{cases}
        x_{h_{f}} & \text{if } f\in (B'\setminus B_{N})\cap W_{a}\\
        x_{g} & \text{if } f= h_{g} \text{ for some } g\in (B'\setminus B_{N})\cap W_{a}\\
        x_{f} & \text{otherwise. }
        \end{cases}
\end{align*}
Also, define $h^{-}\in H_{C}^{-}$ in the following similar way. For every $g\in (B'\setminus B_{N})\cap (F_{\ell}\setminus W_{a})$, choose a $h_{g}\in (F_{\ell}\setminus W_{a})\setminus B'$ such that $h_{g}\neq h_{g'}$ for distinct $g,g'\in (B'\setminus B_{N})\cap (F_{\ell}\setminus W_{a})$. For all $x\in C$, set
\begin{align*}
h^{-}(x)_{f}=
    \begin{cases}
        x_{h_{f}} & \text{if } f\in (B'\setminus B_{N})\cap (F_{\ell}\setminus W_{a})\\
        x_{g} & \text{if } f= h_{g} \text{ for some } g\in (B'\setminus B_{N})\cap (F_{\ell}\setminus W_{a})\\
        x_{f} & \text{otherwise. }
        \end{cases}
\end{align*}
Thus we have $h=h^{+}h^{-}\in H_{C}$, and given any $j,k\in\{1,2,\dots,n\}$,
\begin{align*}
    hC_{j}\cap C_{k}=\{x\mid x_{h_{g}}=i_{g}^{j} \text{ for all } g\in A_{j}\setminus B_{N} \text{ and } x_{f}=i_{f}^{k} \text{ for all } f\in A_{k}\}.
\end{align*}
Thus,
\begin{align*}
    \mu_{C}(hC_{j}\cap C_{k})=\prod_{g\in A_{j}\setminus B_{N}}P_{i_{g}^{j}}\prod_{g\in A_{k}\setminus B_{N}}P_{i_{g}^{k}}
    =\mu_{C}(C_{j})\mu_{C}(C_{k}).
\end{align*} 
Now, let $\epsilon>0$ and let $\Omega$ be any finite collection of measurable subsets of $C$. Note that cylinder sets defined over finitely many elements of $F_{\ell}$ in $C$ generate the $\sigma$-algebra of $C$ inherited from $X_{p}$. This means there is a finite collection $\Omega'$ of cylinder sets defined over finitely many coordinates from $F_{\ell}$ in $C$ such that for every $A\in \Omega$, we can find a $C_{A}\in\Omega'$ such that $\mu_{C}(C_{A}\bigtriangleup A)<\frac{\epsilon}{100}$. So there is an $h\in H_{C}$ such that for any $C_{1},C_{2}\in \Omega'$, we have that 
\begin{align*}
    \mu_{C}(hC_{1}\cap C_{2})=\mu_{C}(C_{1})\mu_{C}(C_{2}).
\end{align*} 
Let $A,B\in \Omega$ and take $C_{A},C_{B}\in \Omega'$ such that $\mu_{C}(C_{A}\bigtriangleup A)<\frac{\epsilon}{100}$ and $\mu_{C}(C_{B}\bigtriangleup B)<\frac{\epsilon}{100}$. Since $\mu_{p}$ is measure preserving, it's not too hard to see that
\begin{align*}
    |\mu_{C}(hA\cap B)-\mu_{C}(A)\mu_{C}(B)|&=|\mu_{C}(hA\cap B)-\mu_{C}(hC_{A}\cap C_{B})\\
    &\hspace*{1cm}+\mu_{C}(C_{A})\mu_{C}(C_{B})-\mu_{C}(A)\mu_{C}(B)|\\ 
    &\leq|\mu_{C}(hA\cap B)-\mu_{C}(hC_{A}\cap C_{B})|\\
    &\hspace*{1cm}+|\mu_{C}(C_{A})\mu_{C}(C_{B})-\mu_{C}(A)\mu_{C}(B)|\\
    &<\frac{\epsilon}{2}+\frac{\epsilon}{2}\\
    &=\epsilon.
\end{align*}
\end{proof} \par
\begin{remark}

Proposition 1.3 implies that there exists a positive integer $M$ and a cylinder set $C\subset X_{p}$ over a ball $B$ in $F_{\ell}$ of radius $M$ centered at the identity element such that 
\begin{align*}
    D=C\cap \{x\mid a^{*}(x)\leq M \text{ and } m^{*}(x)\leq M\}
\end{align*}
has positive measure.
\end{remark}

\begin{proposition} $f$ is constant almost everywhere on $D$.

\end{proposition}
\begin{proof} Choose $\alpha\in\R$ such that $A_{\epsilon}=\{x\mid f(x)-\alpha\mid<\epsilon\}\cap D$ has positive measure for every $\epsilon >0$ and let $B_{\epsilon}=D\setminus A_{\epsilon}$. By way of contradiction, assume that there exists an $\epsilon>0$ such that $\mu_{\textbf{P}}(B_{\epsilon})>0$. By Proposition 3.4, the action of $H_{C}$ on $C$ is weakly mixing and thus ergodic. So there are automorphisms $V^{+}\in H_{C}^{+}$ and $V^{-}\in H_{C}^{-}$ such that, setting $V=V^{+}V^{-}$, we get that $\mu_{p}(VA_{\epsilon}\cap B_{\epsilon})>0$. Set $A=A_{\epsilon}\cap C^{-1}B_{\epsilon}$.
We claim that $J(\mathscr{A}_{p,a},V)=V(\phi^{-1}\mathscr{A}_{q,a},V)=0$. This will finish the proof because, by Proposition 3.3, we get that $f(Vx)=f(x)$ for a.e. $x\in A$ which contradicts that $f$ takes distinct values on $A_{\epsilon}$ and $B_{\epsilon}$.\par
We have $m(x)\leq M$ for all $x\in D$ so that Corollary 1.4 implies that $\phi(V^{-}x)_{g}=\phi(x)_{g}$ for all $g\in W_{a}$ and a.e. $X\in D$. By equation (2) from Proposition 3.2, we have that $J(\mathscr{A}_{p,a},V^{-})=0$ and $J(\phi^{-1}\mathscr{A}_{q,a},V^{-})=J(\mathscr{A}_{q,a},\phi V^{-}\phi^{-1})\circ\phi=0$ a.e. on $A$. Now, let $x'\in A'=V^{-}A$. Since $a_{\phi}\leq M$, we get that $\phi(x')_{g}=\phi(V^{+}x')_{g}$ for $g\in F_{\ell}\setminus W_{a}$ and for a.e. $x\in A'$ by Corollary 1.4. Since $V^{+}$ is measure preserving, we get that $J(\mathscr{A}_{p,a},V)=J(\phi^{-1}\mathscr{A}_{q,a},V)=0$. Thus, the cocycle equation gives
\begin{align*}
    J(\mathscr{A}_{p,a},V)=J(\mathscr{A}_{p,a},V^{+})\circ V^{-}+J(\mathscr{A}_{p,a},V^{-})=0
\end{align*}
and
\begin{align*}
    J(\phi^{-1}\mathscr{A}_{q,a},V)=J(\phi^{-1}\mathscr{A}_{q,a},V^{+})\circ V^{-}+J(\phi^{-1}\mathscr{A}_{q,a},V^{-})=0
\end{align*}
for a.e. $x\in A$.
\end{proof}
\section{Tail Equivalence} 
In the previous section, it was shown that the co-boundary map $f$ is constant on the set $D$, which begets the relation
\begin{align*}J(\mathscr{A}_{p,a},V)=J(\mathscr{A}_{q,a},\phi V\phi^{-1})\circ\phi 
\end{align*}
on $D$ a.e. for every $V\in G_{p}$. In this section, we will show that such a formula holds on a set larger than $D$: the set of all points that are tail equivalent to a point in $D$. Let $f(x)=b$ for a.e. $x\in D$. We set $A=\{x\in X_{p}\mid f(x)=b\}$.
\begin{proposition}
If $V\in H_{C}^{+}\cup H_{C}^{-}$, then $\mu_{p}(VD\setminus A)=0$.
\end{proposition}
\begin{proof}
First, let's suppose that $V\in H_{C}^{-}$. For a.e. $x\in D$, we have that $(Vx)_{g}=x_{g}$ for all $g\in W_{a}\cup B$. Since $m_{\phi}(x)\leq M$, we have that $\phi(V(x))_{g}=\phi(x)_{g}$ for a.e. $x\in D$ and all $g\in W_{a}$ by Corollary 1.4. Equation (2) from Proposition 3.2 shows that
\begin{align*}J(\mathscr{A}_{p,a},V)(x)=J(\mathscr{A}_{q,a},\phi V\phi^{-1})(\phi(x))=0
\end{align*}
a.e. on $D$, and Proposition 3.3 and Proposition 3.6 show that $f(V(x))=f(x)=b$ for a.e. $x\in D$. Thus, $\mu_{p}(VD\setminus A)=0$.\par
Now suppose that $V\in H_{C}^{+}$. Thus we have $(Vx)_{g}=x_{g}$ for all $g\in (F_{\ell}\setminus W_{a})\cup B$ and, since $a_{\phi}(x)\leq M$,  $(\phi(V(x))_{g}=\phi(x)_{g}$ for all $g\in F_{\ell}\setminus W_{a}$. Hence we see that $V$ simply permutes the atoms in the measurable partition coming from $\mathscr{A}_{p,a}$. So we get that
\begin{align*}
    J(\mathscr{A}_{p,a},V)(x)=\log\frac{d\mu_{p}V}{d\mu_{p}}(x)=\log\frac{d\mu_{q}\phi V\phi^{-1}}{d\mu_{q}}(\phi(x))=J(\mathscr{A}_{q,a},\phi V\phi^{-1})(\phi(x))
\end{align*} for a.e. $x\in D$, which also gives that $f(V(x))=f(x)=b$ a.e. on $D$. So the proposition holds.
\end{proof}
\begin{corollary}
There exists a null set $E\subset D$ such that, for every $x\in D\setminus E$, every integer $n\geq 0$, and every $x'\in X_{p}$ with $x_{g}=x_{g}$ for all $g\notin \{a^{k}|-M-n<k<-M\}$ (or for all $g\notin \{a^{k}|M<k<M+n\}$), we have $x'\in A$.
\end{corollary}
\begin{proof} Fix $n\geq 0$. Define subgroups $H_{+}^{n}\subset H_{C}^{+}$ and $H_{-}^{n}\subset H_{C}^{-}$ by 
\begin{align*}
    H_{-}^{n}=\{h\in H_{C}^{-}\mid  (hx)_{g}=x_{g} \text{ for all } g\neq a^{k} \text{ for some } -M-n<k<-M\text{ and for all } x\in X_{p}\}
\end{align*}
and
\begin{align*}
    H_{+}^{n}=\{h\in H_{C}^{+}\mid  (hx)_{g}=x_{g} \text{ for all } g\neq a^{k} \text{ for some } M<k<M+n\text{ and for all } x\in X_{p}\}
\end{align*} for every $x\in X_{p}$. Since these subgroups are finite, Proposition 4.1 implies that there is a null set $E_{n}\subset D$ with $V(x)\in A$ for every $V\in H_{+}^{n} \cup H_{-}^{n}$ and every $x\in D\setminus E_{n}$. The set $E=\cup_{n\geq 1}E_{n}$ will then satisfy the desired property. 
\end{proof}
\section{The Beta Function}
In this section, we will introduce the invariant which will allow us to show that our Bernoulli shifts have the same weights: the beta function. The following definition of the beta function is due to Tuncel \cite{TU} using the definition of pressure given in \cite{WA}. For another formulation of the beta function, see \cite{PT}.
\begin{definition}
Let $(X,d)$ be a compact metric space and $T:X\to X$ be a homeomorphism. Given some $n\in \mathbb{N}$ and $\epsilon>0$, we say a set $E\subset X$ is $(n,\epsilon)\textbf{-separated}$ if whenever $x,y\in E$, where $x\neq y$, there exists some integer $i$ with $0\leq i\leq n-1$ and $d(T^{i}x,T^{i}y)>\epsilon$.  Given a continuous function $f:X\to\mathbb{R}$, we define the \textbf{pressure} $\mathscr{P}_{T}(f)$ of $f$ relative to $T$ by the equation
\begin{align*}
\mathscr{P}_{T}(f)= \lim_{\epsilon\to0}\limsup_{n\to\infty}\frac{1}{n}\log\bigg[\sup_{ E}\sum_{x\in E}e^{\sum_{i=0}^{n-1}f(T^{i}x)}\bigg]
\end{align*}
where $E$ ranges over all $(n,\epsilon)$-separated sets in $X$.
\end{definition}

For continuous functions $f,g:X\to\mathbb{R}$, we have the following:
\begin{enumerate}
\item $\mathscr{P}_{T}(f)=\sup_{\mu}(h(\mu)+\int fd\mu)$
where $\mu$ ranges over all $T$-invariant Borel probability measures on $X$ and $h(\mu)$ is the entropy of $T$ with respect to $\mu$ \cite{TU}. Note that this is one version of the variational principle. 
    \item $\mathscr{P}_{T}(f+g\circ T-g)=\mathscr{P}_{T}(f)$.
\end{enumerate}

\begin{definition}
Let $F_{\ell}\curvearrowright(X_{p},\mu_{p})$ be a Bernoulli shift  with probability vector $p$ over the free group $F_{\ell}$. Fix a generator $a\in F_{\ell}$ and define $T:X_{p}\to X_{p}$ to be the map $T(x)=a^{-1}x$. We define the \textbf{beta function} $\beta_{p,a}:\mathbb{R}\to\mathbb{R}$ by the formula $\beta_{p,a}(t)=\mathscr{P}_{T}(-tI_{\mu_{p}})$, where $I_{\mu_{p}}$ is the conditional information function associated to the past algebra $\mathscr{A}_{p,a}$.
\end{definition}
Combining equation (1) from Proposition 3.2 with Lemma 3 from \cite{PT} yields the following formula for the beta function:
\begin{align*}
   \beta_{p,a}(t)=\lim_{n\to\infty}\bigg(\int\exp((1-t)J(\mathscr{A}_{p,a},T^{n}))d\mu_{p}\bigg)^{\frac{1}{n}}.
\end{align*}

The following theorem appears as Corollary 8 in \cite{PT}.
\begin{theorem}The beta function $\beta_{p,a}$ is analytic, i.e., it is infinitely differentiable everywhere and can be expressed locally about every point as a Taylor series.

\end{theorem}
\begin{proposition}
For a Bernoulli shift $F_{\ell}\curvearrowright(X_{p},\mu_{p})$ defined by a vector $p=[P_{1},\dots,P_{m}]$, we have the following formula for the beta function:
\begin{align*}
    \beta_{p,a}(t)=P_{1}^{t}+\dots + P_{m}^{t}
\end{align*}
for all $t\in\mathbb{R}$.
\end{proposition}
\begin{proof} Let $t\in\mathbb{R}$. Then
\begin{align*}
    \beta_{p,a}(t)
   &=\lim_{n\rightarrow \infty}\bigg(\int\exp((1-t)J(\mathscr{A}_{p,a},T^{n}))d\mu_{p}\bigg)^{1/n} \\
  &\stackrel{(1)}{=} \lim_{n\rightarrow \infty}\bigg(\int\exp\bigg((1-t)\sum_{i=0}^{n-1}-\log(P_{x_{a^{i}}})\bigg)d\mu_{p}(x)\bigg)^{1/n} \\
  &= \lim_{n\rightarrow \infty}\bigg(\int \prod_{i=0}^{n-1} P_{x_{a^{i}}}^{t-1}d\mu_{p}(x)\bigg)^{1/n} \\
  &= \lim_{n\rightarrow \infty}\bigg(\sum_{\sigma\in\{1,\dots,m\}^{\{0,\dots,n-1\}}}((P_{\sigma(0)}\cdots P_{\sigma(n-1)})^{t-1}P_{\sigma(0)}\cdots P_{\sigma(n-1)})\bigg)^{1/n}\\
  &=\lim_{n\rightarrow \infty}\bigg(\sum_{\sigma\in\{1,\dots,m\}^{\{0,\dots,n-1\}}}(P_{\sigma(0)}\cdots P_{\sigma(n-1)})^{t}\bigg)^{1/n}\\
  &=\lim_{n\rightarrow \infty}((P_{1}^{t}+\dots+P_{m}^{t})^{n})^{1/n}\\
   &=P_{1}^{t}+\dots+P_{m}^{t}.
\end{align*}
\end{proof}
The next proposition shows that the formula for the beta function of a Bernoulli shift from Proposition 5.4 gives us an invariant.

\begin{proposition}
Let $p=[P_{1},\dots,P_{m}]$ and $q=[Q_{1},\dots,Q_{n}]$. If \begin{align*}
    P_{1}^{t}+\cdots + P_{m}^{t}=Q_{1}^{t}+\cdots + Q_{n}^{t}
\end{align*}
for all $t\in\mathbb{R}$, then $m=n$ and there is a permutation $\pi \in$ Sym$(1,\dots,n)$ such that $P_{\pi(i)}=Q_{i}$ whenever $1\leq i\leq n$.
\end{proposition}

\begin{proof}
Let $k\in \mathbb{N}$, and set $p_{k}=P_{1}^{k}+\dots + P_{m}^{k}$ and $q_{k}=Q_{1}^{k}+\dots + Q_{n}^{k}$. By hypothesis, we have $p_{k}=q_{k}$ for every $k\in \mathbb{N}$. Let $e_{k,p}$ and $e_{k,q}$ denote the $k$th elementary symmetric polynomials over $p$ and $q$ respectively for every $k\in \mathbb{N}$. By Newton's identities, we have that $p_{k}=q_{k}$ for every $k\in \mathbb{N}$ implies that $e_{k,p}=e_{k,q}$ for all $k\in \mathbb{N}$. Since the polynomials $\prod_{i}^{m}(x-P_{i})$ and $\prod_{i=1}^{n}(x-Q_{i})$ can be expanded as $\sum_{k=0}^{m}(-1)^{k}e_{k,p}x^{n-k}$ and $\sum_{k=0}^{n}(-1)^{k}e_{k,q}x^{n-k}$, we have that $\prod_{i}^{m}(x-P_{i})=\prod_{i=1}^{n}(x-Q_{i})$. This implies that these polynomials have the same roots with the same multiplicities, so that $q$ is just a permutation of the entries of $p$.
\end{proof}
\section{The Main Theorem}
We now arrive at our main theorem. For the proof, we will show, using tail equivalence, that we can restrict the integral in the equation that defines  $\beta_{p,a}$ to a positive measure subset on which the information cocycles for the shifts over their respective past algebras agree up to conjugation by $\phi$. This is enough to show that the beta functions for $p$ and $q$ are equal. First, we recall the following facts about our Bernoulli shifts $F_{\ell}\curvearrowright(X_{p},\mu_{p})$ and $F_{\ell}\curvearrowright(X_{q},\mu_{q})$ given by the vectors $p$ and $q$ respectively. Fix a generator $a\in F_{\ell}$ and suppose $\phi:(X_{p},\mu_{p})\to(X_{q},\mu_{q}) $ is a finitary isomorphism with finite expectation.
Let $C$ be the cylinder set given in Remark 3.5 and $D$ the positive measure subset of $C$ on which $a_{\phi},m_{\phi}\leq M$. Proposition 3.3 gives us a measurable $f:X_{p}\to \mathbb{R}$ that satisfies the co-boundary equation from Proposition 3.3 and is constant on $A$.

\begin{theorem}
Let $F_{\ell}\curvearrowright(X_{p},\mu_{p})$ and $F_{\ell}\curvearrowright(X_{q},\mu_{q})$ be Bernoulli shifts given by two probability vectors $p$ and $q$, and suppose there exists a finitary isomorphism between them with finite expectation. Then $p$ and $q$ are equal up to a permutation of their entries. 
\end{theorem}
For the proof, we will need the following technical lemma.

\begin{lemma}
If $n\geq 1$ and $A=\{x\in X_{p}\mid f(x)=b\}$, then there is a $t_{0}>0$ such that we have
\begin{align*}
    \limsup_{n\to\infty}\bigg(\int_{A\cap a^{n}A}\exp \bigg(tJ(\mathscr{A}_{p,a},a^{-n})\bigg)d\mu_{p}\bigg)^{\frac{1}{n}}
    &=\lim_{n\to\infty}\bigg(\int_{X_{p}}\exp \bigg(tJ(\mathscr{A}_{p,a},a^{-n})\bigg)d\mu_{p}\bigg)^{\frac{1}{n}}\\
    &=\beta_{p,a}(1-t)
\end{align*}
for all $t>t_{0}$.
\end{lemma}
\begin{proof}
 
First, we will partition $X_{p}$ into the equivalence classes induced by the relation of tail equivalence. Indeed, let $M$ be the integer given in Remark 3.6 and choose some $n>2M$. Let $\mathscr{B}$ be the sub-$\sigma$-algebra generated by the collection
\begin{align*}
    \bigcup_{g\in F_{\ell}\setminus\{a^{k}\mid M<k<M+n\}}g^{-1}\alpha_{p}.
\end{align*}
Following Remark 2.7, $\mathscr{B}$ factors $X_{p}$ into a measurable partition $\{[x]_{\mathscr{B}}\}_{x\in X_{p}}$ and a family of conditional measures $\mu_{p}^{\mathscr{B}}$. In fact, each $[x]_{\mathscr{B}}$ is a finite set and each conditional measure $\mu_{p}^{\mathscr{B}}(\cdot|[x]_{\mathscr{B}})$ is an atomic measure such that 
\begin{align*}
    \mu_{p}^{\mathscr{B}}(\{y\}|[x]_{\mathscr{B}})=\prod_{k=M+1}^{M+n-1}P_{y_{a^{k}}}
\end{align*}
for every $y\in[x]_{\mathscr{B}}$. Let $E\subset D$ be the null set given by Corollary 4.2. Note that, if $x\in D\setminus E$, we get that $[x]_{\mathscr{B}}\subset A$. Also, for every $n\geq 1$, if $x\in a^{n}(D\setminus E)$, we get that $a^{-n}[x]_{\mathscr{B}}\subset A$. Hence
\begin{align*}
\mu_{p}^{\mathscr{B}}(A\cap a^{n}A\mid[x]_{\mathscr{B}})=1 \text{ for a.e. }x\in D\cap a^{n}D. \tag{$*$}\end{align*}
We recall that $C$ is a cylinder set such that there exists a ball $B\subset F_{\ell}$ of radius $M$ that is centered at the identity such that $x_{g}=y_{g}$ for all $g\in B$ whenever $x,y\in C$. Let $(i_{g})_{g\in B}\subset \{1,\dots,m\}$ be the tuple over $B$ that determines $C$. Note that $\mu_{p}(D\cap a^{n}D)=\mu_{p}(D)^{2}$ since we chose $n>2M$. This means that we can find a $t_{0}>0$ such that 
\begin{align*}
    \mu_{p}(D\cap a^{n}D)=\mu_{p}(D)^{2}\geq (\prod_{k=-M}^{M}P_{i_{a^{k}}})^{t} \tag{$**$}
\end{align*}
for all $t>t_{0}$. Thus for all $t>t_{0}$:

\begin{align*}
    \int_{A\cap a^{n}A}\exp tJ(\mathscr{A}_{p,a},a^{-n})d\mu_{p}
   & \stackrel{\text {Thm 2.3}}{=} \int  \int_{A\cap a^{n}A}\exp tJ(\mathscr{A}_{p,a},a^{-n})d\mu_{p}^{\mathscr{B}}(\cdot| [x]_{\mathscr{B}})d\mu_{p}(x) \\
    &\stackrel{(*)}{\geq}\int_{D\cap a^{n}D}  \int\exp tJ(\mathscr{A}_{p,a},a^{-n})d\mu_{p}^{\mathscr{B}}(\cdot| [x]_{\mathscr{B}})d\mu_{p}(x) \\
    &=\int_{D\cap a^{n}D}[\sum_{y\in[x]_{\mathscr{B}}}(\prod_{k=M+1}^{n-m-1}P_{y_{a^{k}}})^{1-t}(\prod_{k=-M}^{M}P_{i_{a^{k}}})^{-t}]d\mu_{p}(x)\\
    &= \mu_{p}(D\cap a^{n}D)(\prod_{k=-M}^{M}P_{i_{a^{k}}})^{-t}\int_{x_{p}}\exp\bigg( tJ(\mathscr{A}_{p,a},a^{-(n-2M)})\bigg)d\mu_{p}\\
    &\stackrel{(**)}{\geq}\int_{x_{p}}\exp\bigg( tJ(\mathscr{A}_{p,a},a^{-(n-2M)})\bigg)d\mu_{p}.
\end{align*}

Following from the last display,
\begin{align*}
    \limsup_{n\to\infty}\bigg(\int_{A\cap a^{n}A}\exp\bigg( tJ(\mathscr{A}_{p,a},a^{-n})\bigg)d\mu_{p}\bigg)^{\frac{1}{n}}
    &\geq \limsup_{n\to\infty}\bigg(\int_{x_{p}}\exp \bigg(tJ(\mathscr{A}_{p,a},a^{-(n-2M)})\bigg)d\mu_{p}\bigg)^{\frac{1}{n}}\\
    &=\lim_{n\to\infty}\bigg(\int_{x_{p}}\exp\bigg(tJ(\mathscr{A}_{p,a},a^{-n})\bigg)d\mu_{p}\bigg)^{\frac{1}{n}}\\
    &=\beta_{a,p}(1-t).
    \end{align*}
Thus
\begin{align*}
    \limsup_{n\to\infty} \bigg(\int_{A\cap a^{n}A}\exp\bigg( tJ(\mathscr{A}_{p,a},a^{-n})\bigg)d\mu_{p}\bigg)^{\frac{1}{n}}
   =\beta_{a,p}(1-t)
\end{align*}
holds since the reverse inequality follows from monotonicity of the measure.
\end{proof}\par
Now, we are ready to prove the main theorem.

\begin{proof} It is enough to show that, for every $t\in\mathbb{R}$, $\beta_{p}(t)=\beta_{q}(t)$. Indeed, assuming $\beta_{p}(t)=\beta_{q}(t)$ for all $t\in\mathbb{R}$, we get that $P_{1}^{t}+\dots+P_{m}^{t}=Q_{1}^{t}+\dots+Q_{n}^{t}$ for all $t\in\mathbb{R}$ by Proposition 5.4. This implies that $p$ and $q$ only differ by a permutation using Proposition 5.5. \par
 Proposition 3.3 and Proposition 3.6 imply that 
\begin{align*}
    J(\mathscr{A}_{p,a},a^{-n})(x)=J(\mathscr{A}_{q,a},a^{-n})(\phi(x))
\end{align*}
a.e. on $A\cap a^{-n}A$ for every $n\geq 1$. Furthermore, Lemma 6.2 guarantees a $t_{0}>0$ such that for every $t>t_{0}$
\begin{align*}
    \beta_{p,a}(1-t)
    &=\limsup_{n\to\infty}\bigg(\int_{A\cap a^{n}A}\exp\bigg( tJ(\mathscr{A}_{p,a},a^{-n})\bigg)d\mu_{p}\bigg)^{\frac{1}{n}}\\
    &=\limsup_{n\to\infty}\bigg(\int_{\phi(A)\cap a^{n}\phi(A)}\exp\bigg( tJ(\mathscr{A}_{q,a}),a^{-n})\bigg)d\mu_{q}\bigg)^{\frac{1}{n}}\\
    &\leq \lim_{n\to\infty}\bigg(\int\exp\bigg(tJ(\mathscr{A}_{q,a},a^{-n})\bigg)d\mu_{q}\bigg)^{\frac{1}{n}}=\beta_{q,a}(1-t).
\end{align*} 
By symmetry, we have that $\beta_{p,a}(1-t)=\beta_{q,a}(1-t)$ for all $t>t_{0}$. The analyticity of the $\beta$ function (Theorem 5.3) implies that
$\beta_{p,a}(t)=\beta_{q,a}(t)$ for all $t\in \mathbb{R}$.
\end{proof}
\bibliographystyle{amsplain}
\bibliography{fincodingsbib}
\end{document}